\documentclass[12pt,reqno]{amsart} 
\usepackage{amssymb,amscd,url}
\usepackage{constants} %% \C is next constant, \Cl{name} is labeled constant, \Cr{name} references old constant
\usepackage{xcolor}    %% allows colors, used for edit comments

\begin{document}

% \allowdisplaybreaks

\newcommand{\JOE}[1]{{\textup{\color{blue}~$\bigstar$ \textsf{{\bf Joe:} [#1]}}}}
\newcommand{\WADE}[1]{{\textup{\color{blue}~$\bigstar$ \textsf{{\bf Wade:} [#1]}}}}

%%%%%%%%%%%%%%%%%%%%%%%%%%%%%%%%%%%%%%%%%%%%%%%%%%%%%%%%%%%%%%%%%%%%%%
%% Title and Author Information

\title[The size of semigroup orbits modulo primes]
      {The size of semigroup orbits modulo primes}
\date{\today}

\author{Wade Hindes}
\email{wmh33@txstate.edu}
\address{Department of Mathematics, Texas State University,  
		601 University Dr., 
		San Marcos, TX 78666 USA}

\author[Joseph H. Silverman]{Joseph H. Silverman}
\email{joseph\_silverman@brown.edu}
\address{Mathematics Department, Box 1917
         Brown University, Providence, RI 02912 USA}
%% \author{CoAuthor}
%% \email{CoAuthor}
%% \address{CoAuthor}

\thanks{
  Silverman's research supported by Simons Collaboration Grant \#712332, and by National Science Foundation Grant \#1440140 while in residence at MSRI in spring 2023.
}

%%%%%%%%%%%%%%%%%%%%%%%%%%%%%%%%%%%%%%%%%%%%%%%%%%%%%%%%%%%%%%%%%%%%%%

\hyphenation{ca-non-i-cal semi-abel-ian}

%%%%%%%%%%%%%%%%%%%%%%%%%%%%%%%%%%%%%%%%%%%%%%%%%%%%%%%%%%%%%%%%%%%%%%
% Theorem environments

\newtheorem{theorem}{Theorem}
\newtheorem{lemma}[theorem]{Lemma}
\newtheorem{sublemma}[theorem]{Sublemma}
\newtheorem{conjecture}[theorem]{Conjecture}
\newtheorem{proposition}[theorem]{Proposition}
\newtheorem{corollary}[theorem]{Corollary}
\newtheorem*{claim}{Claim}

\theoremstyle{definition}
% The * surpresses numbering
\newtheorem{definition}[theorem]{Definition}
\newtheorem{example}[theorem]{Example}
\newtheorem{remark}[theorem]{Remark}
\newtheorem{question}[theorem]{Question}

\theoremstyle{remark}
\newtheorem*{acknowledgement}{Acknowledgements}

%%%%%%%%%%%%%%%%%%%%%%%%%%%%%%%%%%%%%%%%%%%%%%%%%%%%%%%%%%%%%%%%%%%%%%

%%%%%%%% Set Up Environment for Notation %%%%%%%%%%%%%%
% This is currently set to allow quite wide items to be defined
\newenvironment{notation}[0]{%
  \begin{list}%
    {}%
    {\setlength{\itemindent}{0pt}
     \setlength{\labelwidth}{4\parindent}
     \setlength{\labelsep}{\parindent}
     \setlength{\leftmargin}{5\parindent}
     \setlength{\itemsep}{0pt}
     }%
   }%
  {\end{list}}

%%%%%%%% Set Up Environment for Parts in Theorems %%%%%%%%%%%%%%
\newenvironment{parts}[0]{%
  \begin{list}{}%
    {\setlength{\itemindent}{0pt}
     \setlength{\labelwidth}{1.5\parindent}
     \setlength{\labelsep}{.5\parindent}
     \setlength{\leftmargin}{2\parindent}
     \setlength{\itemsep}{0pt}
     }%
   }%
  {\end{list}}
% Use \Part{(a)}, instead of \item[(a)], to ensure upright font
\newcommand{\Part}[1]{\item[\upshape#1]}

%%%%%%%% Set Up Macro for Cases %%%%%%%%%%%%%%
\def\Case#1#2{%
\paragraph{\textbf{\boldmath Case #1: #2.}}\hfil\break\ignorespaces}

%%%%%%%%%%%%%%%%%%
% Greek Alphabet %
%%%%%%%%%%%%%%%%%%
\renewcommand{\a}{\alpha}
\renewcommand{\b}{\beta}
\newcommand{\g}{\gamma}
\renewcommand{\d}{\delta}
\newcommand{\e}{\epsilon}
\newcommand{\f}{\varphi}
\newcommand{\bfphi}{{\boldsymbol{\f}}}
\renewcommand{\l}{\lambda}
\renewcommand{\k}{\kappa}
\newcommand{\lhat}{\hat\lambda}
\newcommand{\m}{\mu}
\newcommand{\bfmu}{{\boldsymbol{\mu}}}
\renewcommand{\o}{\omega}
\renewcommand{\r}{\rho}
\newcommand{\rbar}{{\bar\rho}}
\newcommand{\s}{\sigma}
\newcommand{\sbar}{{\bar\sigma}}
\renewcommand{\t}{\tau}
\newcommand{\z}{\zeta}

\newcommand{\D}{\Delta}
\newcommand{\G}{\Gamma}
\newcommand{\F}{\Phi}
\renewcommand{\L}{\Lambda}

%%%%%%%%%%%%%%%%%%%%
% Fraktur Alphabet %
%%%%%%%%%%%%%%%%%%%%
\newcommand{\ga}{{\mathfrak{a}}}
\newcommand{\gb}{{\mathfrak{b}}}
\newcommand{\gn}{{\mathfrak{n}}}
\newcommand{\gp}{{\mathfrak{p}}}
\newcommand{\gq}{{\mathfrak{q}}}
\newcommand{\gA}{{\mathfrak{A}}}
\newcommand{\gB}{{\mathfrak{B}}}
\newcommand{\gC}{{\mathfrak{C}}}
\newcommand{\gD}{{\mathfrak{D}}}
\newcommand{\gP}{{\mathfrak{P}}}
\newcommand{\gS}{{\mathfrak{S}}}

%%%%%%%%%%%%%%%%%%%
% Barred Alphabet %
%%%%%%%%%%%%%%%%%%%
\newcommand{\Abar}{{\bar A}}
\newcommand{\Ebar}{{\bar E}}
\newcommand{\kbar}{{\bar k}}
\newcommand{\Kbar}{{\bar K}}
\newcommand{\Pbar}{{\bar P}}
\newcommand{\Sbar}{{\bar S}}
\newcommand{\Tbar}{{\bar T}}

%%%%%%%%%%%%%%%%%%%%%%%%%
% Calligraphic Alphabet %
%%%%%%%%%%%%%%%%%%%%%%%%%
\newcommand{\Acal}{{\mathcal A}}
\newcommand{\Bcal}{{\mathcal B}}
\newcommand{\Ccal}{{\mathcal C}}
\newcommand{\Dcal}{{\mathcal D}}
\newcommand{\Ecal}{{\mathcal E}}
\newcommand{\Fcal}{{\mathcal F}}
\newcommand{\Gcal}{{\mathcal G}}
\newcommand{\Hcal}{{\mathcal H}}
\newcommand{\Ical}{{\mathcal I}}
\newcommand{\Jcal}{{\mathcal J}}
\newcommand{\Kcal}{{\mathcal K}}
\newcommand{\Lcal}{{\mathcal L}}
\newcommand{\Mcal}{{\mathcal M}}
\newcommand{\Ncal}{{\mathcal N}}
\newcommand{\Ocal}{{\mathcal O}}
\newcommand{\Pcal}{{\mathcal P}}
\newcommand{\Qcal}{{\mathcal Q}}
\newcommand{\Rcal}{{\mathcal R}}
\newcommand{\Scal}{{\mathcal S}}
\newcommand{\Tcal}{{\mathcal T}}
\newcommand{\Ucal}{{\mathcal U}}
\newcommand{\Vcal}{{\mathcal V}}
\newcommand{\Wcal}{{\mathcal W}}
\newcommand{\Xcal}{{\mathcal X}}
\newcommand{\Ycal}{{\mathcal Y}}
\newcommand{\Zcal}{{\mathcal Z}}

%%%%%%%%%%%%%%%%%%%%%%%%%%%%
% Blackboard Bold Alphabet %
%%%%%%%%%%%%%%%%%%%%%%%%%%%%
\renewcommand{\AA}{\mathbb{A}}
\newcommand{\BB}{\mathbb{B}}
\newcommand{\CC}{\mathbb{C}}
\newcommand{\FF}{\mathbb{F}}
\newcommand{\GG}{\mathbb{G}}
\newcommand{\NN}{\mathbb{N}}
\newcommand{\PP}{\mathbb{P}}
\newcommand{\QQ}{\mathbb{Q}}
\newcommand{\RR}{\mathbb{R}}
\newcommand{\ZZ}{\mathbb{Z}}

%%%%%%%%%%%%%%%%%%%%%%%%%%
% Boldface Math Alphabet %
%%%%%%%%%%%%%%%%%%%%%%%%%%
\newcommand{\bfa}{{\boldsymbol a}}
\newcommand{\bfb}{{\boldsymbol b}}
\newcommand{\bfc}{{\boldsymbol c}}
\newcommand{\bfd}{{\boldsymbol d}}
\newcommand{\bfe}{{\boldsymbol e}}
\newcommand{\bff}{{\boldsymbol f}}
\newcommand{\bfg}{{\boldsymbol g}}
\newcommand{\bfi}{{\boldsymbol i}}
\newcommand{\bfj}{{\boldsymbol j}}
\newcommand{\bfp}{{\boldsymbol p}}
\newcommand{\bfr}{{\boldsymbol r}}
\newcommand{\bfs}{{\boldsymbol s}}
\newcommand{\bft}{{\boldsymbol t}}
\newcommand{\bfu}{{\boldsymbol u}}
\newcommand{\bfv}{{\boldsymbol v}}
\newcommand{\bfw}{{\boldsymbol w}}
\newcommand{\bfx}{{\boldsymbol x}}
\newcommand{\bfy}{{\boldsymbol y}}
\newcommand{\bfz}{{\boldsymbol z}}
\newcommand{\bfA}{{\boldsymbol A}}
\newcommand{\bfF}{{\boldsymbol F}}
\newcommand{\bfB}{{\boldsymbol B}}
\newcommand{\bfD}{{\boldsymbol D}}
\newcommand{\bfG}{{\boldsymbol G}}
\newcommand{\bfI}{{\boldsymbol I}}
\newcommand{\bfM}{{\boldsymbol M}}
\newcommand{\bfP}{{\boldsymbol P}}
\newcommand{\bfzero}{{\boldsymbol{0}}}
\newcommand{\bfone}{{\boldsymbol{1}}}

%%%%%%%%%%%%%%%%%%%%%%%%%%%%%%
% Miscellaneous New Commands %
%%%%%%%%%%%%%%%%%%%%%%%%%%%%%%
\newcommand{\Aut}{\operatorname{Aut}}
\newcommand{\codim}{\operatorname{codim}}
\newcommand{\Crit}{\operatorname{Crit}}
\newcommand{\CritVal}{\operatorname{CritVal}}
\newcommand{\Curve}{{\Gamma}}
\newcommand{\Disc}{\operatorname{Disc}}
\newcommand{\Div}{\operatorname{Div}}
\newcommand{\Dom}{\operatorname{Dom}}
\newcommand{\End}{\operatorname{End}}
\newcommand{\Fbar}{{\bar{F}}}
\newcommand{\Fix}{\operatorname{Fix}}
\newcommand{\Gal}{\operatorname{Gal}}
\newcommand{\GL}{\operatorname{GL}}
\newcommand{\Hom}{\operatorname{Hom}}
\newcommand{\Image}{\operatorname{Image}}
\newcommand{\Isom}{\operatorname{Isom}}
\newcommand{\hhat}{{\hat h}}
\newcommand{\Ker}{{\operatorname{ker}}}
\newcommand{\limstar}{\lim\nolimits^*}
\newcommand{\limstarn}{\lim_{\hidewidth n\to\infty\hidewidth}{\!}^*{\,}}
\newcommand{\Mat}{\operatorname{Mat}}
\newcommand{\maxplus}{\operatornamewithlimits{\textup{max}^{\scriptscriptstyle+}}}
\newcommand{\MOD}[1]{~(\textup{mod}~#1)}
\newcommand{\Mor}{\operatorname{Mor}}
\newcommand{\Moduli}{\mathcal{M}}
\newcommand{\Norm}{{\operatorname{\mathsf{N}}}}
\newcommand{\notdivide}{\nmid}
\newcommand{\normalsubgroup}{\triangleleft}
\newcommand{\NS}{\operatorname{NS}}
\newcommand{\onto}{\twoheadrightarrow}
\newcommand{\ord}{\operatorname{ord}}
\newcommand{\Orbit}{\mathcal{O}}
\newcommand{\Orb}{\operatorname{Orb}}
\newcommand{\Per}{\operatorname{Per}}
\newcommand{\PrePer}{\operatorname{PrePer}}
\newcommand{\PGL}{\operatorname{PGL}}
\newcommand{\Pic}{\operatorname{Pic}}
\newcommand{\Prob}{\operatorname{Prob}}
\newcommand{\Proj}{\operatorname{Proj}}
\newcommand{\Qbar}{{\overline{\QQ}}}
\newcommand{\rank}{\operatorname{rank}}
\newcommand{\Rat}{\operatorname{Rat}}
\newcommand{\Resultant}{\operatorname{Res}}
\renewcommand{\setminus}{\smallsetminus}
\newcommand{\sgn}{\operatorname{sgn}} 
\newcommand{\SL}{\operatorname{SL}}
\newcommand{\Span}{\operatorname{Span}}
\newcommand{\Spec}{\operatorname{Spec}}
\renewcommand{\ss}{\textup{ss}}
\newcommand{\stab}{\textup{stab}}
\newcommand{\Support}{\operatorname{Supp}}
\newcommand{\tors}{{\textup{tors}}}
\newcommand{\tr}{{\textup{tr}}} 
\newcommand{\Trace}{\operatorname{Trace}}
\newcommand{\trianglebin}{\mathbin{\triangle}} % symmetric set difference
\newcommand{\UHP}{{\mathfrak{h}}}    % Upper half plane
\newcommand{\<}{\langle}
\renewcommand{\>}{\rangle}

\newcommand{\pmodintext}[1]{~\textup{(mod}~#1\textup{)}}
\newcommand{\ds}{\displaystyle}
\newcommand{\longhookrightarrow}{\lhook\joinrel\longrightarrow}
\newcommand{\longonto}{\relbar\joinrel\twoheadrightarrow}
\newcommand{\SmallMatrix}[1]{%
  \left(\begin{smallmatrix} #1 \end{smallmatrix}\right)}

%%%%%%%%%%%%%%%%%%%%%%%%%%%%%%%%%%%%%%%%%%%%%%%%%%%%%%%%%%%%%%%%%%%%%%

\begin{abstract} 
Let $V$ be a projective variety defined over a number field $K$, let $S$ be a polarized set of endomorphisms of $V$ all defined over~$K$, and let $P\in V(K)$. For each prime $\mathfrak{p}$ of $K$, let $m_{\mathfrak{p}}(S,P)$ denote the number of points in the orbit of $P\bmod\mathfrak{p}$ for the semigroup of maps generated by $S$. Under suitable hypotheses on $S$ and $P$, we prove an analytic estimate for $m_{\mathfrak{p}}(S,P)$ and use it to show that the set of primes for which $m_{\mathfrak{p}}(S,P)$ grows subexponentially as a function of $\operatorname{\mathsf{N}}_{K/\mathbb{Q}}\mathfrak{p}$ is a set of density zero. For $V=\PP^1$ we show that this holds for a generic set of maps $S$ provided that at least two of the maps in $S$ have degree at least four.
\end{abstract}

\maketitle

\tableofcontents

%%%%%%%%%%%%%%%%%%%%%%%%%%%%%%%%%%%%%%%%%%%%%%%%%%
\section{Introduction}
%%%%%%%%%%%%%%%%%%%%%%%%%%%%%%%%%%%%%%%%%%%%%%%%%%
A general heuristic in arithmetic dynamics over number fields is that the dynamical systems generated by ``unrelated'' self-maps~$f_1,f_2: V\rightarrow V$ should not be too similar. For example, they should not have identical canonical heights~\cite{MR2352570}, they should not have infinitely many common preperiodic points~\cite{MR2817647,MR4088354,MR3335283}, their orbits should not have infinite intersection~\cite{MR2922378}, and arithmetically their orbits should not have unexpectedly large common divisors~\cite{MR4050112}. It is not always clear what  ``unrelated'' should mean, but in any case it includes the assumption that~$f_1$ and~$f_2$ do not share a common iterate.
\par
Similarly, we expect that the points in semigroup orbits generated by all finite compositions of ``unrelated'' maps~$f_1$ and~$f_2$ should be asymptotically large~\cite{BellHindesNote,MR4449715} when ordered by height, where now unrelated means that the semigroup is not unexpectedly small. For example, the semigroup is small if it contains no free sub-semigroups requiring at least~$2$ generators; cf.~\cite{Tits}. 
\par
In this note, we study the size of semigroup orbits over finite fields. In particular, we show that under suitable hypotheses that a free semigroup of maps defined over a number field generates many large orbits when reduced modulo primes. See~\cite{MR2549537,MR3319121,MR3767353,MR2448661} for additional results in this vein.  

% The paper~\cite{MR3767353} of Chang, D'Andrea, Ostafe, Shparlinski, Sombra, studies a parameterized family of multi-variable polynomials, but it deals with the lengths of the orbits of the individual polynomials, not the full orbit of the monoid that they generate.

\begin{definition}
\label{definition:KQVS}
We set notation that will remain in effect throughout this note.
\begin{center}
\begin{tabular}{cl}
$K/\QQ$ & a number field\\
$V/K$ & a smooth projective variety defined ove r$K$ \\
$r\ge1$ & an integer \\
$S=\{f_1,\dots,f_r\}$ & a set of morphisms~$f_i:V\rightarrow V$ defined over $K$ \\
$d_1,\ldots,d_r$ & real numbers satisfying $d_i>1$ \\
$\Lcal\in\Pic_K(V)\otimes\RR$ & line bundles satisfying $f_i^*\Lcal\cong\Lcal^{\otimes d_i}$\\
$M_S$ & the semigroup generated by~$S$ under composition\\
$\Orb_S(P)$ & the orbit $\{f(P): f\in M_S\}$ of a point $P\in V$ \\
\end{tabular}
\end{center}
\end{definition}

The following property will play a crucial role in some of our results. 

\begin{definition} 
\label{definition:wanderperiodic}
A point~$P\in V$ is called \emph{strongly $S$-wandering} if the evaluation map
\begin{equation}
\label{eqn:wandering}
M_S\longrightarrow V,\quad f\longmapsto f(P),
\end{equation}
is injective. 
\end{definition} 

\begin{remark} 
If~$V=\PP^1$ and~$S$ is any sufficiently generic set of maps as described in Section~\ref{sec:examples}, then the set of points that fail to be strongly $S$-wandering is a set of bounded height. In particular, it follows in this case that all infinite orbits contain strongly wandering points, and this weaker condition is sufficient for our orbit bounds. 
\end{remark}

Our goal is to study the number of points in the reduction of~$\Orb_S(P)$ modulo primes. We set some additional notation, briefly recall a standard definition,
and then define our principal object of study.

\begin{center}
\begin{tabular}{cl}
$R_K$ & the ring of integers of~$K$\\
$\Spec(R_K)$ & the set of prime ideals of $R_K$ \\
$\Norm\gp$ & ${}:=\#R_K/\gp$, the norm of $\gp\in\Spec(R_K)$ \\
%% $m_\gp$ & ${}= m_\gp(S,P) := \#\Orb_{\tilde S}(\tilde P\bmod\gp)$ \\
\end{tabular}
\end{center}

\begin{definition}
Let $\gp\in\Spec(R_K)$. A $K$-morphism $f:V\to{V}$ has \emph{good reduction at~$\gp$} if there is a scheme $\Vcal_\gp/R_\gp$ that is proper and smooth over~$R_\gp$ and an $R_\gp$-morphism $F_\gp:\Vcal_\gp\to\Vcal_\gp$ whose generic fiber is $f:V\to{V}$. We note that in this situation, reduction modulo~$\gp$ of~$f$ commutes with iteration of~$f$.
\end{definition}

\begin{definition}
Let $\gp\in\Spec(R_K)$. If all of the maps~$f_1,\dots,f_r$ have good reduction modulo~$\gp$, then we write
\[
m_\gp := m_\gp(S,P) = \#\Orb_{\tilde S}(\tilde P\bmod\gp)
\]
for the size of~$\Orb_S(P)$ when it is reduced mod~$\gp$. Otherwise we formally set~$m_\gp=\infty$. 
\end{definition}
 
Our main result is an analytic formula that implies that~$m_\gp$ is not too small on average.

\begin{theorem}
\label{theorem:sumlogpoverpmpe}
Assume that~$M_S$ is a free semigroup, that~$P\in{V}(K)$ is a strongly $S$-wandering point, and that~$r=\#S\ge2$. 
Then there exists a constant~$\Cl{t1}=\Cr{t1}(K,V,S,P)$ such that for all~$\e>0$,  
\begin{equation}
\label{eqn:mainthmsum}
\sum_{\gp\in\Spec R_K}
\frac{\log\Norm\gp}{\Norm\gp\cdot m_\gp(S,P)^\e} 
\le \Cr{t1} \e^{-1}.
\end{equation}
\end{theorem}

\begin{remark}
The principal result of the paper~\cite{MR2448661} is an estimate exponentially weaker than~\eqref{eqn:mainthmsum} in the case that~$r=\#S=1$, while a principal result of the paper~\cite{MR1395936} is an estimate that exactly mirrors~\eqref{eqn:mainthmsum} with~$m_\gp$ equal to the number of points on the  mod~$\gp$ reduction of the multiples of a point on an abelian variety. Thus the present paper, as well as the papers~\cite{BellHindesNote,MR4449715}, suggest that the analogy proposed
\[
\left(\begin{tabular}{@{}l@{}}
arithmetic of points\\
of an abelian variety\\
\end{tabular}\right)
\quad\Longleftrightarrow\quad
\left(\begin{tabular}{@{}l@{}}
arithmetic of points in orbits\\
of a dynamical system\\
\end{tabular}\right)
\]
described in~\cite{MR4007163} and~\cite[\S6.5]{MR2884382} may be more accurate when the dynamical system on the right-hand side is generated by at least two non-commuting maps, rather than using orbits coming from iteration of a single map.
\end{remark}

The estimate~\eqref{eqn:mainthmsum} can be used to show that there are few primes~$\gp$ for which~$m_{\gp}$ is subexponential compared to~$\Norm\gp$. We quantify this assertion in the following corollary.

\begin{corollary}
\label{corollary:densitymplep} 
Let~$S$,~$M_S$ and~$P$ be as in Theorem~\textup{\ref{theorem:sumlogpoverpmpe}}, and let~$\overline\d$ and~$\d$ denote the \textup(upper\textup) logarithmic analytic density on sets of primes as described in Definition~\textup{\ref{definition:analyticdensity}}.%
\begin{parts}
\Part{(a)}
There is a constant~$\Cl{2}=\Cr{2}(K,V,S,P)$ such that
\[
\overline\d\Bigl(\bigl\{\gp\in\Spec R_K : m_\gp\le\Norm\gp^\g \bigr\}\Bigr)
\le \Cr{2}\g
\]
holds for all~$0<\g<1$. \vspace{.15cm}
\Part{(b)}
Let~$L(t)$ be a subexponential function, i.e., a function with the property
that
\[
\lim_{t\to\infty} \frac{L(t)}{t^\mu}=0
\quad\text{for all~$\mu>0$.}
\]
Then
\[
\d\Bigl(\bigl\{\gp\in\Spec R_K : m_\gp\le L(\Norm\gp) \bigr\}\Bigr) = 0.
\]
\end{parts}
\end{corollary}

In the special case that~$V=\PP^1$, we show that the conclusions of Theorem~\ref{theorem:sumlogpoverpmpe} and Corollary~\ref{corollary:densitymplep} are true for generic sets of maps. In the statement of the next result, we write~$\Rat_d$ for the space of rational maps of~$\PP^1$ of degree~$d\geq2$, so in particular~$\Rat_d$ is an affine variety of dimension~$2d+1$; see~\cite[\S4.3]{SilvDyn} for details. 

\begin{theorem}
\label{thm:generic} 
Let $r\ge2$, and let $d_1,\ldots,d_r$ be integers satisfying 
\[
d_1,d_2\geq4 \quad\text{and}\quad d_3,\ldots,d_r\ge2.
\]
Then there is a Zariski dense subset
\[
\Ucal=\Ucal(d_1,\dots,d_r)\subseteq\Rat_{d_1}\times\dots\times \Rat_{d_r}
\]
such that the inequality~\eqref{eqn:mainthmsum} in 
Theorem~\textup{\ref{theorem:sumlogpoverpmpe}} and the density estimates in Corollary~\textup{\ref{corollary:densitymplep}} are true for all number fields~$K/\QQ$, all~$S\in\Ucal(K)$, and all~$P\in\PP^1(K)$ for which~$\Orb_S(P)$ is infinite. 
\end{theorem}

The contents of this paper are as follows. In Section~\ref{sec:orbitsmodp} we build upon prior work~\cite{MR1395936,MR2448661} of the second author to prove Theorem~\ref{theorem:sumlogpoverpmpe} and Corollary~\ref{corollary:densitymplep}. Then in Section~\ref{sec:examples} we use results from~\cite{MR2922378,MR4449715,MR4449722} to construct many sets of maps on~$\PP^1$ for which the bounds in Section~\ref{sec:orbitsmodp} apply. The key step is to construct a point in every infinite orbit that is strongly wandering. The construction is explicit, and in particular, Theorem~\ref{thm:ratmaps} describes an explicit set~$\mathcal{U}$ for which Theorem~\ref{thm:generic} is true.

%%%%%%%%%%%%%%%%%%%%%%%%%%%%%%%%%%%%%%%%%%%%%%%%%%%%%%%%%%%%%%%%%%%%%%
\section{The Size of Orbits Modulo~$\gp$}
\label{sec:orbitsmodp}
%%%%%%%%%%%%%%%%%%%%%%%%%%%%%%%%%%%%%%%%%%%%%%%%%%%%%%%%%%%%%%%%%%%%%%
We start with a key estimate. 

\begin{proposition}
\label{proposition:loglogDmlelogm}
For each~$m\ge2$, we define an integral ideal
\begin{equation}
  \label{eqn:gDmdef}
  \gD(m) = \gD(m;K,V,S,P) := \prod_{\substack{\gp\in\Spec R_K\\m_\gp(S,P)\le m\\}} \gp.
\end{equation}
There are constants~$C_i=C_i(K,V,S,P)$ such that the following hold\textup:%
\begin{parts}
\Part{(a)}
If~$r=\#S=1$, then
\[
\log\log\Norm\gD(m)\le \Cl{m1} m
\quad\text{for all~$m\ge2$.}
\]
\Part{(b)}
Assume that~$S$ generates a free semigroup,
that~$P\in{V}(K)$ is strongly $S$-wandering,
and that~$r=\#S\ge2$. Then
\[
\log\log\Norm\gD(m)\le \Cl{m2} \log m
\quad\text{for all~$m\ge2$.}
\]
\end{parts}
\end{proposition}

\begin{proof}
(a)\enspace
This is~\cite[Proposition~10]{MR2448661}.
\par\noindent(b)\enspace  
For notational convenience, we use the standard combinatorics
notation~$[r]=\{1,2,\ldots,r\}$.  Also to ease notation, 
we write
\[
f_{\bfi} := f_{i_1}\circ f_{i_2}\circ\cdots\circ f_{i_k}
\quad\text{for}\quad \bfi=(i_1,\ldots,i_k)\in[r]^k.
\]
\par
Let
\[
m \ge 1 \quad\text{and}\quad
k=k(m):=\left\lceil\dfrac{\log(m+1)}{\log r}\right\rceil.
\]
For each good reduction prime~$\gp$, we consider the map
that sends a function~$f_\bfi$ to the image of~$P$ under
reduction modulo~$\gp$,
\begin{equation}
\label{eqn:redmodpmaprktoOPp}
  [r]^k\longrightarrow \Orbit_{\tilde S}(\tilde P\bmod\gp),\quad
  \bfi \longmapsto f_\bfi(\tilde P\bmod\gp).
\end{equation}
If
\[
m_\gp\le m,\quad\text{then}\quad r^k > m_\gp
\quad\text{by our choice of~$k$,}
\]
so the map~\eqref{eqn:redmodpmaprktoOPp} cannot be injective
(pigeonhole principle) and there exist
\[
\text{$\bfi\ne\bfj$ in~$[r]^k$}
\quad\text{satisfying}\quad
f_\bfi(\tilde P\bmod\gp) = f_\bfj(\tilde P\bmod\gp).
\]
Sine we have assumed that~$P$ is strongly wandering, i.e., that the map
\[
M_S \longrightarrow V(K),\quad
f\longmapsto f(P),
\]
is injective, it follows that the global points are distinct,
\[
f_\bfi(P)\ne f_\bfj(P),
%% f_{i_1}\circ f_{i_2}\circ\cdots\circ f_{i_k}(P)
%% \ne
%% f_{j_1}\circ f_{j_2}\circ\cdots\circ f_{j_k}(P), 
\]
so the ideals generated by their ``differences'' are non-zero.
\par
More formally,
let~$\gC=\gC_K$ be the ideal of~$R_K$ that depends only
on~$K$ and is described in~\cite[Lemma~9]{MR2448661}.
Then for~$P\in{X(K)}$, we can write
\[
f_\bfi(P) = \bigl[A_0(\bfi),\ldots,A_N(\bfi)\bigr]
\]
with~$A_0(\bfi),\ldots,A_N(\bfi)\in{R_K}$ and such that the ideal
\[
\gA(\bfi) := A_0(\bfi)R_K+\cdots+A_N(\bfi)R_K
\quad\text{divides the ideal~$\gC$.}
\]
Then for~$\gp\nmid\gC$ we have
\begin{align*}
  f_\bfi(P)&\equiv f_\bfj(P)  \pmodintext{\gp} \\
  &\Longleftrightarrow\quad
  A_u(\bfi)A_v(\bfj) \equiv A_v(\bfi)A_u(\bfj) \pmodintext{\gp}
  \quad\text{for all~$0\le u,v\le N$.}  
\end{align*}
We define a difference ideal
\[
  \gB(\bfi,\bfj) := \sum_{0\le u,v\le N}
  \Bigl( A_u(\bfi)A_v(\bfj) - A_v(\bfi)A_u(\bfj) \Bigr) R_K, \\
\]
and the product of the difference ideals
\[
  \gD'(m) :=   \prod_{\bfi,\bfj\in[r]^k,\;\bfi\ne\bfj} \gB(\bfi,\bfj).
\]
Then
\[
\gp\nmid\gC~\text{and}~m_\gp\le m
\quad\Longrightarrow\quad
\gp\mid \gD'(m),
\]
and hence
\[
\gD(m)\mid\gC\cdot\gD'(m).
\]
Since~$\gC$ depends only on~$K$, it
remains to estimate the norm of~$\gD'(m)$.
\par
Using~\cite[Proposition~7]{MR2448661}, we get
\[
\frac{1}{[K:\QQ]}\log \frac{\Norm\gB(\bfi,\bfj)}{\Norm\gA(\bfi)\cdot\Norm\gA(\bfj)}
\le
h\bigl( f_\bfi(P) \bigr) + h\bigl( f_\bfj(P) \bigr) + \Cl{0}.
\]
Since~$\Norm\gA(\bfi)$ and~$\Norm\gA(\bfj)$ are smaller
than~$\Norm\gC$, we find that
\[
\frac{1}{[K:\QQ]}\log \Norm\gB(\bfi,\bfj)
\le
h\bigl( f_\bfi(P) \bigr) + h\bigl( f_\bfj(P) \bigr) + \Cl{2z}
\]
Next we apply the height estimate
\[
h\bigl( f_\bfi(P) \bigr) \le \Cl{1} \cdot \prod_{u=1}^k d_{i_u},
\]
which is a weak form of~\cite[Lemma~2.1]{MR3958050}.
This yields
\[
\frac{1}{[K:\QQ]}\log \Norm\gB(\bfi,\bfj)
\le
\Cr{1} \cdot \prod_{u=1}^k d_{i_u} + \Cr{1} \cdot \prod_{u=1}^k d_{j_u} + \Cr{2z}.
\]  
This gives
\begin{align*}
\log \gD'(m)
&= \sum_{\bfi,\bfj\in[r]^k,\;\bfi\ne\bfj} \log \gB(\bfi,\bfj) \\
&\le \sum_{\bfi,\bfj\in[r]^k,\;\bfi\ne\bfj}
\Bigl( \Cr{1}\cdot \prod_{u=1}^k d_{i_u}
+ \Cr{1}\cdot \prod_{u=1}^k d_{j_u} + \Cr{2z} \Bigr) \\
&\le \Cl{3}\cdot r^k \cdot \sum_{ \bfi \in [r]^k } \prod_{u=1}^k d_{i_u} \\
&= \Cr{3}\cdot \Bigl( r\cdot \sum_{i\in[r]} d_i \Bigr)^k \\
&\le \Cl{4} \cdot \biggl( r\cdot \sum_{i\in[r]} d_i \biggr)^{1+\frac{\log(m+1)}{\log r}}.
\end{align*}
Hence
\[
\log\log \gD'(m) \le \Cl{5} \cdot \log(m+1) + \Cl{6}.
\]
Since~$m\ge2$, we can absorb~$\Cr{6}$ into~$\Cr{5}$, although we remark
that if we leave in~$\Cr{6}(K,V,S,P)$, then we can
take~$\Cr{5}$ to depend on only the degrees of the maps in~$S$,
\[
\Cr{5} = \Cr{5}(d_1,\ldots,d_r) =
1 + \frac{\log(d_1+\cdots+d_r)}{\log r}.
\]
This completes the proof of Proposition~\ref{proposition:loglogDmlelogm}.
\end{proof}

\begin{proof}[Proof of Theorem~$\ref{theorem:sumlogpoverpmpe}$]
To ease notation, we let
\[
g(t) = \frac{\log t}{t} \quad\text{and}\quad G(t)=\frac{1}{t^\e}.
\]
We start with two elementary estimates. First, the Mean Value Theorem
gives
\begin{equation}
  \label{eqn:mvtestimate}
  G(m)-G(m+1)
  \le \sup_{m\le t\le m+1} -G'(t)
  = \sup_{m\le t\le m+1} \frac{\e}{t^{1+\e}} = \frac{\e}{m^{1+\e}}.
\end{equation}
Second, an easy integral calculation gives
\begin{equation}
  \label{eqn:intbypartsestimate}
  \sum_{m\ge1} g(m)G(m) \le \int_1^\infty \frac{\log x}{x^{1+\e}} = \frac{1}{\e^2}.
\end{equation}
We use these and our other calculations to estimate
\begin{align*}
  \sum_{\gp\in\Spec R_K} & \frac{\log\Norm\gp}{\Norm\gp\cdot m_\gp^\e} \\
  &=  \sum_{\gp\in\Spec R_K} g(\Norm\gp) \cdot G(m_\gp)
  \quad\text{by definition of~$g$ and~$G$,}\\
  &= \sum_{m\ge1} G(m) \sum_{\substack{\gp\in\Spec R_K\\ m_\gp=m\\}} g(\Norm\gp) \\
  &= \sum_{m\ge1} \Bigl(G(m)-G(m+1)\Bigr)
  \sum_{\substack{\gp\in\Spec R_K\\ m_\gp\le m\\}} g(\Norm\gp)
  \quad\text{Abel summation,} \\
  &\le \sum_{m\ge1} \frac{\e}{m^{1+\e}}
  \sum_{\substack{\gp\in\Spec R_K\\ m_\gp\le m\\}} g(\Norm\gp)
  \quad\text{from \eqref{eqn:mvtestimate},} \\
  &= \sum_{m\ge1} \frac{\e}{m^{1+\e}}
  \sum_{\substack{\gp\in\Spec R_K\\ \gp\mid\gD(m)\\}} g(\Norm\gp)
  \quad\text{by definition \eqref{eqn:gDmdef} of~$\gD(m)$,} \\
  &\le \sum_{m\ge1} \frac{\e}{m^{1+\e}}
  \cdot \Bigl( \Cl{7}\log\log\gD(m)+\Cl{8} \Bigr)
  \quad\text{from~\cite[Corollary~2.3]{MR1395936},} \\
  &\le \Cl{9} \sum_{m\ge1} \frac{\e}{m^{1+\e}}\cdot\log m
  \quad\text{from Proposition~\ref{proposition:loglogDmlelogm}(b),} \\
  &= \Cr{9}\cdot \e\cdot \sum_{m\ge1} g(m)\cdot G(m)
  \quad\text{by definition of~$g$ and~$G$,} \\
  &\le \Cl{10}\e^{-1}
  \quad\text{from \eqref{eqn:intbypartsestimate}.}  
\end{align*}
This completes the proof of Theorem~\ref{theorem:sumlogpoverpmpe}.
\end{proof}

\begin{definition}
\label{definition:analyticdensity}
Let $\Pcal\subset\Spec R_K$ be a set of primes. The \emph{upper logarithmic analytic density of~$\Pcal$} is
\[
\overline\d(\Pcal):=
\limsup_{s\to1^+} (s-1)\sum_{\gp\in\Pcal} \frac{\log\Norm\gp}{\Norm\gp^s}
\]
Similarly, the \emph{logarithmic analytic density of~$\Pcal$}, denoted~~$\d(P)$, is given by the same formula with a limit, instead of a limsup.
\end{definition}

\begin{proof}[Proof of Corollary~$\ref{corollary:densitymplep}$]
(a)\enspace
For any~$0<\g<1$, we let
\[
\Pcal_\g := \bigl\{\gp\in\Spec R_K : m_\gp\le\Norm\gp^\g \bigr\}.
\]
Then
\begin{align}
  \label{eqn:Cegesump1ge}
  \frac{\Cr{t1}}{\e}
  &\ge \sum_{\gp\in\Spec R_K} \frac{\log\Norm\gp}{\Norm\gp\cdot m_\gp(S,P)^\e}
  \quad\text{from Theorem~\ref{theorem:sumlogpoverpmpe},} \notag\\
  &\ge \sum_{\gp\in\Pcal_\g} \frac{\log\Norm\gp}{\Norm\gp\cdot m_\gp(S,P)^\e}
  \quad\text{summing over a smaller set,} \notag\\
  &\ge \sum_{\gp\in\Pcal_\g} \frac{\log\Norm\gp}{\Norm\gp^{1+\g\e}}
  \quad\text{by definition of~$\Pcal_\g$.}
\end{align}
This allows us to estimate the upper logarithmic density of~$\Pcal_\g$ by
\begin{align*}
\overline\d(\Pcal_\g)
&= \limsup_{s\to1^+} (s-1)\sum_{\gp\in\Pcal_\g} \frac{\log\Norm\gp}{\Norm\gp^s} \\
&= \limsup_{\e\to0^+} \g\e \sum_{\gp\in\Pcal_\g} \frac{\log\Norm\gp}{\Norm\gp^{1+\g\e}}
\quad\text{setting~$s=1+\g\e$,} \\
&\le  \limsup_{\e\to0^+} \g\e \cdot \frac{\Cr{t1}}{\e}
\quad\text{from \eqref{eqn:Cegesump1ge},} \\
&= \Cr{t1}\g.  
\end{align*}
This completes the proof of Corollary~\ref{corollary:densitymplep}(a).
\par\noindent(b)\enspace
We let
\[
\Pcal_L := \bigl\{\gp\in\Spec R_K : m_\gp\le L(\Norm\gp) \bigr\}.
\]
The assumption that~$L$ is subexponential means that for all~$\mu>0$
there exists a constant~$\Cl{mu}(L,\mu)$ depending only on~$L$ and~$\mu$
such that
\[
L(t) \le t^\mu \quad\text{for all~$t>\Cr{mu}(L,\mu)$.}
\]
We also note that
\begin{align}
  \label{eqn:mplepmypgtCmu}
  \gp\in\Pcal_L
  &\quad\Longleftrightarrow\quad
  m_\gp \le L(\Norm\gp) \notag\\
  &\quad\Longrightarrow\quad
  m_\gp \le (\Norm\gp)^\mu\quad\text{for all~$\Norm\gp>\Cr{mu}(L,\mu)$.}
\end{align}
We now fix a~$\mu>0$ and estimate
\begin{align*}
  \overline\d(\Pcal_L)
  &= \limsup_{\l\to0^+} \l\sum_{\gp\in\Pcal_L} \frac{\log\Norm\gp}{\Norm\gp^{1+\l}} \\
  &= \limsup_{\l\to0^+}
  \l\sum_{\hidewidth\substack{\gp\in\Pcal_L\\\Norm\gp\ge\Cr{mu}(L,\mu)\\}\hidewidth}
  \frac{\log\Norm\gp}{\Norm\gp^{1+\l}} \\
  &\omit\hspace*{3em}\hfill
  \text{since~$\mu$ is fixed, so we can discard finitely many terms,} \\
  &\le   \limsup_{\l\to0^+}
  \l\sum_{\hidewidth\substack{\gp\in\Pcal_L\\\Norm\gp\ge\Cr{mu}(L,\mu)\\}\hidewidth}
  \frac{\log\Norm\gp}{\Norm\gp} \cdot \frac{1}{m_\gp^{\l/\mu}}
  \quad\text{from \eqref{eqn:mplepmypgtCmu},} \\  
  &\le   \limsup_{\l\to0^+}
  \l\sum_{\gp\in\Spec R_K}
  \frac{\log\Norm\gp}{\Norm\gp} \cdot \frac{1}{m_\gp^{\l/\mu}} \\
  &\le   \limsup_{\l\to0^+}
  \l \cdot \Cr{t1}\cdot\left(\frac{\l}{\mu}\right)^{-1}
  \quad\text{from Theorem~\ref{theorem:sumlogpoverpmpe},} \\
  &= \Cr{t1}\mu.
\end{align*}
This estimate holds for all~$\mu>0$, so we find that
\[
\overline\d(\Pcal_L) \le \inf_{\mu>0} \Cr{t1}\cdot\mu = 0,
\]
which completes the proof that~$\d(\Pcal_L)=0$. 
\end{proof}

%%%%%%%%%%%%%%%%%%%%%%%%%%%%%%%%%%%%%%%%%
\section{Orbits of Generic Families of Maps of $\PP^1$}
\label{sec:examples} 
%%%%%%%%%%%%%%%%%%%%%%%%%%%%%%%%%%%%%%%%%
In this section, we show that there are many sets of endomorphisms of~$\PP^1$ for which Theorem~\ref{theorem:sumlogpoverpmpe} holds. To make this statement precise, we need some definitions.

\begin{definition}
Let~$f$ be a non-constant rational map of~$\PP^1$ defined over~$\Qbar$. A point~$w\in\PP^1(\Qbar)$ is a \emph{critical value} of~$f$ if~$f^{-1}(w)$ contains fewer than~$\deg(f)$ elements. It is a \emph{simple critical value} if
\[
\#f^{-1}(w) = \deg(f)-1.
\]
The map~$f$ is \emph{critically simple} if all of its critical values are simple
\end{definition}
  
\begin{definition}
\label{def:crit}
Let~$f$ and~$g$ be non-constant rational maps of~$\PP^1$ with respective critical value sets~$\CritVal_{f}$ and~$\CritVal_{g}$. We say that~$f$ and~$g$ are \emph{critically separated} if 
\[
\CritVal_{f}\cap\CritVal_{g}=\varnothing.
\]
\end{definition}

Our first result says that the conclusions of Theorem~\ref{theorem:sumlogpoverpmpe} and Corollary~\ref{corollary:densitymplep} hold for certain sets~$S$ that contain a pair of critically simple and critically separated maps and initial points~$P$ with infinite orbit.  

\begin{theorem}
\label{thm:ratmaps} 
Let~$K/\QQ$ be a number field, let~$S$ be a set of endomorphisms of~$\PP^1$ defined over~$K$ containing a pair of critically simple and critically separated maps of degree at least~$4$, and let~$P\in\PP^1(K)$ be a point with infinite~$S$-orbit. Then there is a constant~$\Cl{3z}=\Cr{3z}(K,S,P)$ such that for all~$\epsilon>0$, 
\[
\sum_{\gp\in\Spec_{R_K}}\frac{\log\Norm\gp}{\Norm\gp\cdot m_\gp(S,P)^\e}\leq \Cr{3z}\cdot\epsilon^{-1}. 
\]
\end{theorem}

\begin{remark} 
In particular, there is a constant~$\Cl{hPS}(S)$ such that Theorem~\ref{thm:ratmaps} holds for all~$P\in\PP^1(K)$ satisfying~$h(P)>\Cr{hPS}(S)$; see Lemma~\ref{lem:hts}. 
\end{remark} 

We start with a definition and some basic height estimates.

\begin{definition} 
\label{definition:moderateS}
A point $P\in V$ is \emph{moderately $S$-preperiodic} if
\[
g\circ f(P)=f(P)\quad\text{for some $f,g\in M_S$ with $g\ne1$.}
\]
\end{definition}

\begin{lemma}
\label{lem:hts}
Let $V/\Qbar$ be a variety, and let~$S=\{f_1,\dots, f_r\}$ be a set of polarized endomorphisms as described in Definition~\textup{\ref{definition:KQVS}}. Then there exists a constant~$\Cl{SVL}=\Cr{SVL}(S,V,\Lcal)$ such that the following statements hold for all~$Q\in{V}(\Qbar)$\textup: 
\begin{parts} 
\Part{(a)}
If~$Q$ is moderately $S$-preperiodic
as described in Definition~\textup{\ref{definition:moderateS}}, then
\[
h_{\Lcal}(Q)\leq \Cr{SVL}.
\]
In particular, this is true if~$\Orb_S(Q)$ is finite.
\Part{(b)}
If~$h_{\Lcal}(Q)>\Cr{SVL}$, then
\[
h_{\Lcal}\bigl(f(Q)\bigr)\geq h_{\Lcal}(Q)
\quad\text{for all~$f\in M_S$.} 
\]
\end{parts}    
\end{lemma} 
\begin{proof}
These estimates are proven in~\cite[Lemma~2.10]{BellHindesNote}.  
\end{proof}

We combine Lemma~\ref{lem:hts} with the techniques in~\cite{MR4449715,MR2774592} to obtain the following result for pairs of maps that are critically simple and critically separated.

\begin{proposition}
\label{prop:wanderingpts} 
Let~$f_1$ and~$f_2$ be endomorphisms of~$\PP^1$ of degree at least~$4$, let~$S=\{f_1,f_2\}$, and suppose that~$f_1$ and~$f_2$ are critically simple and critically separated. 
\begin{parts} 
\Part{(a)}
The semigroup~$M_{S}$ is free.  
\Part{(b)}
Let~$P\in\PP^1(\Qbar)$ be a point whose~$S$-orbit~$\Orb_S(P)$ is infinite. Then there exists a point~$Q\in\Orb_S(P)$ such that~$Q$ is strongly $S$-wandering as described in Definition~\textup{\ref{definition:wanderperiodic}}.
\end{parts}  
\end{proposition}
\begin{proof} 
(a)\enspace
See~\cite[Proposition 4.1]{MR4449715}. 
\par\noindent(b)\enspace
We fix a number field~$K$ over which~$P$,~$f_1$, and~$f_2$ are defined. Letting $\D\subset\PP^1\times\PP^1$ be the diagonal, we define three curves
\begin{align*}
\Curve_{i}&:=(f_i\times f_i)^{-1}(\Delta)\quad\text{for $i=1,2$},\\
\Curve_{1,2}&:=(f_1\times f_2)^{-1}(\Delta).
\end{align*}
Then the main results in~\cite{MR2774592} (see also~\cite[Proposition~4.6]{MR4449715}) imply that the curves~$\Curve_{1}$ and~$\Curve_{2}$ are each the union of~$\Delta$ and an irreducible curve of geometric genus~$\ge2$, while~$\Curve_{1,2}$ is itself an irreducible curve of geometric genus~$\ge 2$.
\par
More specifically, the assumption that~$f_1$ and~$f_2$ are critically simple implies from~\cite[Corollary 3.6]{MR2774592} that \text{$\mathcal{C}_{1}\setminus\Delta$}  and \text{$\mathcal{C}_{2}\setminus\Delta$} are irreducible, while the assumption that~$f_1$ and~$f_2$ are critically separated implies from \cite[Proposition 3.1]{MR2774592} that $\mathcal{C}_{1,2}$ is irreducible. It then follows from~\cite[pages~208 and~210]{MR2774592} that the geometric genera of these curves are given by the formulas
\begin{align*}
\operatorname{genus}(\Curve_{i}\setminus\D) 
&= \bigl(\deg(f_i)-2)\bigr)^2\quad\text{for $i=1,2$}, \\
\operatorname{genus}(\Curve_{1,2}\setminus\D) 
&=(\deg(f_1)-1)(\deg(f_2)-1).
\end{align*}
In particular, the assumption that~$f_1$ and~$f_2$ have degree at least~$4$ ensures that these genera are at least~$2$.
\par
We now invoke Faltings's theorem~\cite{MR718935},~\cite[Theorem~E.0.1]{MR1745599} to deduce that the set
\[
\Sigma:=\Curve_{1,2}(K)
\cup(\Curve_1\setminus\Delta)(K)
\cup(\Curve_2\setminus\Delta)(K)
\] 
is finite. We note that the definition of~$\Sigma$ says that for all~$P,Q\in\PP^1(K)$, we have
\begin{equation}
\label{eqn:fiPeqfjQ}
\left[
\begin{aligned}
P\ne Q\quad\text{and}\quad f_1(P)=f_1(Q)
&\quad\Longrightarrow\quad (P,Q)\in\Sigma \\
P\ne Q\quad\text{and}\quad f_2(P)=f_2(Q)
&\quad\Longrightarrow\quad (P,Q)\in\Sigma \\
f_1(P)=f_2(Q)
&\quad\Longrightarrow\quad (P,Q)\in\Sigma 
\end{aligned}
\right].
\end{equation}
\par
Let $\pi_1,\pi_2:\PP^1\times\PP^1\to\PP^1$ be the two projection maps, and let
\[
\Cl{hLSigma} := 
\max\left(\Bigl\{ h_\Lcal(P) : P\in\pi_1(\Sigma) \Bigr\}
\cup
\Bigl\{ h_\Lcal(P) : P\in\pi_1(\Sigma) \Bigr\}\right)
\]
be the maximum of the heights of the coordinates of the finitely many points in~$\Sigma$. We then set
\[
\Cl{BCCprime} := \max\{ \Cr{SVL},\Cr{hLSigma} \},
\]
where~$\Cr{SVL}$  is the constant that appears in Lemma~\ref{lem:hts}. 
\par
The fact that $\Orb_S(P)\subseteq\PP^1(K)$ is infinite, combined with Northcott's theorem~\cite{MR34607} saying that~$\PP^1(K)$ has only finitely many points of bounded height, implies that there exists a point~$Q\in\Orb_S(P)$ satisfying
\[
h_{\Lcal}(Q)>\Cr{BCCprime}.
\]
\par
We claim that~$Q$ is strongly wandering for~$S$. To see this, 
suppose that
\begin{equation}
\label{eqn:fifjnmQ}
f_{i_1}\circ\dots\circ f_{i_n}(Q)
= f_{j_1}\circ\dots\circ f_{j_m}(Q),
\end{equation}
where without loss of generality we may assume that~$n\ge{m}$.
Our goal is to prove that~$m=n$ and~$i_k=j_k$ for all~$1\le{i}\le{n}$.
\par
To ease notation, we let
\begin{equation}
\label{eqn:fifjnmQFG}
F=f_{i_2}\circ\dots\circ f_{i_n}
\quad\text{and}\quad
G=f_{j_2}\circ\dots\circ f_{j_m}
\end{equation}
be the compositions with the initial map omitted. Thus~\eqref{eqn:fifjnmQ} and~\eqref{eqn:fifjnmQFG} say that
\begin{equation}
\label{eqn:fi1FQfj1GQ}
f_{i_1}\bigl(F(Q)\bigr)= f_{j_1}\bigl(G(Q)\bigr).
\end{equation}
It follows from~\eqref{eqn:fi1FQfj1GQ} and~\eqref{eqn:fiPeqfjQ} that one of the following is true: 
\begin{parts}
\Part{(1)}
$i_1=j_1$ and $F(Q)=G(Q)$.
\Part{(2)}
$i_1=j_1$ and $F(Q)\ne G(Q)$ and $\bigl(F(Q),G(Q)\bigr)\in\Sigma$.
\Part{(3)}
$i_1\ne j_1$  and $\bigl(F(Q),G(Q)\bigr)\in\Sigma$.
\end{parts}
On the other hand, we know that
\[
\bigl(F(Q),G(Q)\bigr)\in\Sigma
\quad\Longrightarrow\quad
h_\Lcal\bigl(F(Q)\bigr)\leq \Cr{BCCprime} <h_\Lcal(Q),
\]
which contradicts Lemma~\ref{lem:hts}. Hence~(2) and~(3) are false, so~(1) is true. 
\par
We recall that~$m\le{n}$, so repeating this argument, we conclude that
\[
i_k=j_k\quad\text{for all $1\le k\le m$.}
\]
If~$m<n$ is a strict inequality, then we obtain
\[
f_{i_{m+1}}\circ\dots\circ f_{i_n}(Q)=Q,
\]
so Lemma~\ref{lem:hts} implies that~$h_\Lcal(Q)\le\Cr{SVL}\le\Cr{BCCprime}$, which is a contradiction. Hence~$m=n$ and~$i_k=j_k$ for all~$1\leq k\leq n$, which completes the proof that~$Q$ is a strongly $S$-wandering point.
\end{proof}

We now have the tools in place to prove Theorem~\ref{thm:ratmaps}. 

\begin{proof}[Proof of Theorem~\textup{\ref{thm:ratmaps}}] 
Let~$S$ be the given set of endomorphisms of~$\PP^1$, and let~$f_1$ and~$f_2$ be the given maps in~$S$ that have degree at least~$4$ and that are critically simple and critically separated. We let
\[
S'=\{f_1,f_2\}.
\]
We are given that the point~$P\in\PP^1(K)$ has infinite~$S$-orbit, and hence by Northcott's theorem~\cite{MR34607}, there are points of arbitrarily large height in~$\Orb_S(P)$. We choose a point
\[
Q'\in \Orb_S(P) \quad\text{satisfying}\quad
h_\Lcal(Q') > \Cr{SVL}(S'),
\]
where~$\Cr{SVL}(S')$ is the constant  associated to the set~$S'$ appearing in Lemma~\ref{lem:hts}.  In particular, it follows from Lemma~\ref{lem:hts}(b) and Northcott's theorem that $\Orb_{S'}(Q')$ must be infinite. Then Proposition~\ref{prop:wanderingpts} implies that~$M_{S'}$ is free and that there is a point
\[
Q\in\Orb_{S'}(Q')\subseteq\Orb_S(P)
\]
that is strongly $S'$-wandering. Applying Theorem~\ref{theorem:sumlogpoverpmpe} to the set~$S'$ and the point~$Q$, we deduce that
\[
\sum_{\gp\in\Spec_{R_K}}\frac{\log\Norm\gp}{\Norm\gp\cdot m_\gp(S,P)^\e}\leq\sum_{\gp\in\Spec_{R_K}}\frac{\log\Norm\gp}{\Norm\gp\cdot m_{\gp}(S',Q)^{\epsilon}}\leq \Cr{3z}\epsilon^{-1} 
\]
for some constant~$\Cr{3z}$ depending on~$S$,~$Q$ (and so~$P$) and~$K$. For this last conclusion, we have also used the fact that
\[
m_{\gp}(S',Q)\leq m_{\gp}(S,P),
\]
which is immediate from the inclusion~$\Orb_{S'}(Q)\subseteq\Orb_S(P)$.     
\end{proof} 

\begin{proof}[Proof of Theorem~\textup{\ref{thm:generic}}]
We recall that~$\Rat_d$ denotes the space of rational maps of degree~$d$. It is shown in~\cite{MR2774592} that if~$d_1,d_2\ge4$, then the set
\[
\Vcal_{d_1,d_2} :=
\left\{ (f_1,f_2)\in \Rat_{d_1}\times\Rat_{d_2} : 
\begin{tabular}{@{}l@{}}
$f_1$ and $f_2$ are critically simple\\
 and critically separated\\
\end{tabular}
\right\}
\]
is Zariski dense in~$\Rat_{d_1}\times\Rat_{d_2}$. Then for any~$d_3,\ldots,d_r\ge2$, the set
\[
\Ucal_{d_1,\ldots,d_r} := \Vcal_{d_1,d_2}
\times \Rat_{d_3} \times\cdots\times \Rat_{d_r}
\]
is Zariski dense in~$\Rat_{d_1}\times\cdots\times\Rat_{d_r}$, and Theorem~\ref{thm:ratmaps} gives us that the desired inequality~\eqref{eqn:mainthmsum} for every~$S$ generated by a set of maps
\[
(f_1,\ldots,f_r)\in\Ucal_{d_1,\ldots,d_r}.
\]
\end{proof}

We conclude with a variant of Theorem~\ref{thm:ratmaps} in which the maps are polynomials. We start with a definition.

\begin{definition}
\label{definition:grouplike}
A polynomial~$f(x)\in\Qbar[x]$ is \emph{power-like} if there exist polynomials~$R(x),C(x),L(x)\in\Qbar[x]$ such that
\begin{align*}
f &= R\circ C\circ L, & \deg(L) &= 1, & \deg(C) &\ge 2, \\
C(x) &= \text{a power map or a Chebyshev polynomial.} \hidewidth
\end{align*}
\end{definition}

\begin{theorem}
\label{thm:poly} 
Let~$K/\QQ$ be a number field, let~$S$ be a set of endomorphisms of~$\PP^1$ defined over~$K$, and let~$P\in\PP^1(K)$ be a point such that~$\Orb_S(P)$ is infinite. Suppose further that~$S$ contains polynomials~$f_1(x),f_2(x)\in{K}[x]$ having the following properties\textup:
\begin{parts} 
\Part{(1)}
Neither~$f_1$ nor~$f_2$ is power-like; see Definition~\textup{\ref{definition:grouplike}}.
\Part{(2)}
For all $g\in\Qbar[x]$ satisfying $\deg(g)\ge2$, we have
\[
f_1\ne f_2\circ g \quad\text{and}\quad f_2\ne f_1\circ g.
\]
\end{parts}
Then there is a constant~$\Cl{5z}=\Cr{5z}(K,S,P)$ such that for all~$\epsilon>0$,
\begin{equation*}
\sum_{\gp\in\Spec_{R_K}}\frac{\log\Norm\gp}{\Norm\gp\cdot m_\gp(S,P)^\e}\leq \Cr{5z}\cdot\epsilon^{-1}. 
\end{equation*} 
\end{theorem}

The proof of Theorem~\ref{thm:poly} is similar to the proof of Theorem~\ref{thm:ratmaps}, except that we use~\cite{MR2922378,MR4449722}  instead of~\cite{MR4449715,MR2774592}. As a first step, we need the following result, which is a polynomial analogue of Proposition~\ref{prop:wanderingpts}.   

\begin{proposition}
\label{prop:polywanderingpts} 
Let~$f_1$ and~$f_2$ be polynomials satisfying the hypotheses of Theorem~\textup{\ref{thm:poly}}, and let~$S=\{f_1,f_2\}$. 
\begin{parts} 
\Part{(a)} The semigroup~$M_{S}$ is free. \vspace{.15cm} 
\Part{(b)} Let~$P\in\PP^1(\Qbar)$ be a point whose $S$-orbit~$\Orb_S(P)$ is infinite. Then there exists a strongly $S$-wandering point~$Q\in\Orb_S(P)$.
\end{parts}
\end{proposition}
\begin{proof} 
(a)\enspace 
See~\cite[Proposition 4.5]{MR4449715}. 
\par\noindent(b)\enspace
The proof is very similar to the proof of Proposition~\ref{prop:wanderingpts}, so we just give a brief sketch, highlighting the differences. We note that we have picked a coordinate function~$x$ on~$\PP^1$. We let~$\infty\in\PP^1$ be the pole of~$x$ and let~$\AA^1=\PP^1\setminus\{\infty\}$. Replacing~$P$ with another point in~$\Orb_S(P)$ if necessary, we may assume that~$P\ne\infty$ is not the point at infinity. We choose a set~$\gS$ of primes of~$K$ so that the ring of~$\gS$-integers~$R_{K,\gS}$ satisfies
\[
P\in \AA^1(R_{K,\gS}) \quad\text{and}\quad
f_1(x),f_2(x)\in R_{K,\gS}[x].
\]
\par
We use the map $f_1\times{f_2}:\AA^2\to\AA^2$ to define three affine curves,
\[
\Curve_{1}:=(f_1\times f_1)^{-1}(\Delta),\quad
\Curve_{2}:=(f_2\times f_2)^{-1}(\Delta),\quad
\Curve_{1,2}:=(f_1\times f_2)^{-1}(\Delta).
\]
Then~\cite[Proposition 4.5]{MR4449715}, itself a consequence of the main results of~\cite{MR2922378,MR4449722}, tell us that these are geometrically irreducible curves of geometric genus at least~$1$. (This is where we use the assumptions~(1) and~(2) on~$f_1$ and~$f_2$.) The Siegel--Mahler theorem for integral points on affine curves~\cite[Theorem~D.9.1]{MR1745599} then implies that
\[
\Curve_{1}(R_{K,\gS}),\quad
\Curve_{2}(R_{K,\gS}),\quad\text{and}\quad
\Curve_{1,2}(R_{K,\gS})\quad\text{are finite sets,}
\]
and hence that
\[
\Sigma:=\Curve_{1,2}(R_{K,\gS})\cup(\Curve_1\setminus\Delta)(R_{K,\gS})\cup(\Curve_2\setminus\Delta)(R_{K,\gS})
\] 
is finite. 
\par
The remainder of the proof of Proposition~\ref{prop:polywanderingpts} is identical to the proof of Proposition~\ref{prop:wanderingpts}, starting with the three possibilities described in~\eqref{eqn:fiPeqfjQ}.
\end{proof}

\begin{proof}[Proof of Theorem~\textup{\ref{thm:poly}}]
The proof of Theorem~\ref{thm:poly} is identical to the proof of Theorem~\ref{thm:ratmaps}. We first use Lemma~\ref{lem:hts}, Proposition~\ref{prop:polywanderingpts}, and the fact that~$\Orb_S(P)$ is infinite to find a point~$Q\in\Orb_S(P)$ that is strongly wandering for~$S'=\{f_1,f_2\}$. We then apply Theorem~\ref{theorem:sumlogpoverpmpe} to the point~$Q$ and the set~$S'=\{f_1,f_2\}$ to deduce the desired result for~$P$ and~$S$. 
\end{proof}

%% \begin{acknowledgement}
%% The author would like to thank ***
%% for their helpful advice.
%% \end{acknowledgement}

%% \bibliographystyle{plain}
%% \bibliography{OrbitsModP}

\end{document}